\documentclass[a4paper,12pt]{amsart}
\usepackage[dvips]{graphicx,color}
\usepackage{amsmath}
\usepackage{amssymb}
\usepackage{amsbsy}
\usepackage{amsgen}
\usepackage{amsopn}
\usepackage{amscd}
\usepackage{amstext}
\usepackage{amsthm}

\theoremstyle{definition}
\newtheorem{dfn}{Definition}[section]
\theoremstyle{plain}
\newtheorem{thm}{Theorem}[section]
\newtheorem*{mthm}{Main Theorem}

\theoremstyle{remark}
\newtheorem{rem}{Remark}[section]

\theoremstyle{definition}

\begin{document}

\title{A bicomplex of Khovanov homology for colored Jones polynomial}
\author{Noboru Ito}\thanks{The author was a Research Fellow of the Japan Society for the Promotion of Science.  This work was partly supported by KAKENHI.  This work was supported in part by IRTG 1529}

\begin{abstract}
We construct a bicomplex for the categorification of the colored Jones polynomial.  This work is motivated by the problem suggested by Anna Beliakova and Stephan Wehrli who discussed the categorification of the colored Jones polynomial in their paper.  
\end{abstract}
\maketitle

\section{Introduction.  }
In the papers \cite{khovanov}, \cite{mt} and \cite{bw}, the categorification of the colored Jones polynomial are given by using the cabling formula (\ref{eq-cj}) in \cite[Theorem 4.15]{km}.  
They introduced the colored Jones polynomial $J_{\bf{n}}$ by \cite[Section 3.1]{bw} and \cite[Section 1.2]{khovanov} ; the notion of a {\it{cabling}} and the cabling formula of colored Jones polynomial are first introduced by \cite{murakami}).  

Let ${\bf{n}}$ $=$ $(n_{1}, \dots, n_{l})$ be the tuple of non-negative integers.  
For an arbitrary oriented framed $l$ component link $L$ whose $i$-th component is colored by the $(n_{i}+1)$-dimensional irreducible representation of $\mathcal{U}_{q}({sl}_{2})$.  Let $J(L^{\bf{n}})$ be the Jones polynomial of ${\bf{n}}$-cable of $L$.  Forming the $m$-cable of a component of $L$, we orient the strands by alternating the original and opposite directions.  The colored Jones polynomial $J_{\bf{n}}$ is given by 
\begin{equation}\label{eq-cj}
J_{\bf{n}}(L) = \sum_{{\bf{k}} = {\bf{0}}}^{\lfloor \frac{{\bf{n}}}{2} \rfloor}(-1)^{|\bf{k}|}\begin{pmatrix} {\bf{n}}-{\bf{k}} \\ {\bf{k}}\end{pmatrix} J(L^{\bf{n} - \bf{2k}})
\end{equation}
where $|{\bf{k}}|$ $=$ $\sum_{i}k_{i}$ and 
\begin{equation}
\begin{pmatrix} {\bf{n}}-{\bf{k}} \\ {\bf{k}}\end{pmatrix} = \prod_{i = 1}^{l} \begin{pmatrix} n_{i} - k_{i} \\ k_{i} \end{pmatrix}.  
\end{equation}
Let us consider the Khovanov complex of a link diagram $D$ where it is defined in \cite[Page 1215, Definition 3]{jacobsson} and denoted by $\mathcal{C}^{i, j}(D)$ as in \cite[Page 237]{viro} or \cite{ito3}.  Set $\mathcal{C}^{i}(D)$ $:=$ $\oplus_{j}\mathcal{C}^{i, j}(D)$.  
For a diagram $D$ of a framed link $L$, the diagram of ${\bf{n}}$-cable $L^{\bf{n}}$ of the fixed $D$ is denoted by $D^{\bf{n}}$.  Following \cite[Section 6.1.1]{w} or \cite[Section 3.2]{bw}, we consider the graph $\Gamma_{\bf{n}}$ corresponding to $D^{\bf{n}}$.  The binomial coefficient $\left(\begin{smallmatrix} {\bf{n}} - {\bf{k}} \\ {\bf{k}} \end{smallmatrix}\right)$ is the number of ways to select ${\bf{k}}$ pairs of neighbourhoods on $l$ lines.  We call such a selection of ${\bf{k}}$ pairs a {\bf{k}}-{\it pairing} for the sake of simplicity.  For a given {\bf{k}}-pairing ${\bf{s}}$, $D^{\bf{s}}$ denotes the cable diagram consisting of components corresponding to unpaired dots.  By definition, $D^{\bf{s}}$ is isotopic to $D^{{\bf{n}} - 2{\bf{k}}}$; see \cite[Page 62]{w}.  Let $I_{\bf{k}}$ be the set of all {\bf{k}}-pairings \cite[Page 66]{w}, $\mathcal{C}_{\bf{n}}^{k, i, j}(D)$ $:=$ $\oplus_{{\bf{s}} \in I_{\bf{k}}, |{\bf{k}}| = k}\mathcal{C}^{i, j}(D^{\bf{s}})$ and $\mathcal{C}_{\bf{n}}^{k, i}(D)$ $:=$ $\oplus_{{\bf{s}} \in I_{\bf{k}}, |{\bf{k}}| = k}$$\mathcal{C}^{i}(D^{\bf{s}})$.  
By using (\ref{eq-cj}), it holds that
\begin{equation}
J_{\bf{n}}(L) = \sum_{k}(-1)^{k}\sum_{\begin{smallmatrix} {\bf{s}} \in I_{\bf{k}}\\ |{\bf{k}}| = k \end{smallmatrix}}\sum_{i, j}(-1)^{i}q^{j}\,\mathrm{rk}\,\mathcal{C}^{i, j}(D^{{\bf{n}} - 2{\bf{k}}}).  
\end{equation}
Then, we have
\begin{equation}
J_{\bf{n}}(L) = \sum_{k}(-1)^{k+i}q^{j}\,\mathrm{rk}\,\mathcal{C}^{k, i, j}_{\bf{n}}(D).  
\end{equation}
The grade $k$ becomes one of the homological grade in the following sense.  
The graph $\Gamma_{\bf{n}}$ corresponding to $D^{\bf{n}}$ is denoted by $\Gamma_{\bf{n}}(D)$.  Let $F^{k}(\Gamma(D))$ be the free Abelian group generated by 
$\left\{ \Gamma_{\bf{s}}(D)~\text{corresponding to}~{\bf{s}} \in I_{\bf{k}}~\text{with}~|{\bf{k}}| = k \right\}$.  

As in \cite[Section 3.3, third paragraph]{bw}, with an edge of $\Gamma_{\bf{n}}$ connecting a ${\bf{k}}$-pairing {\bf{s}} and a ${\bf{k'}}$-pairing ${\bf{s'}}$ we associate a homomorphism $e$ : $F^{k}(\Gamma_{\bf{n}}(D))$ $\to$ $F^{k+1}(\Gamma_{\bf{n}}(D))$ given by gluing an annulus between the strands of the cable which form a pair in ${\bf{s'}}$ but not in ${\bf{s}}$.  
We define the differential $d_{\bf{n}}^{k} :$ $F^{k}(\Gamma_{\bf{n}}(D))$ $\to$ $F^{k+1}(\Gamma_{\bf{n}}(D))$ by $(-1)^{({\bf{s}}, {\bf{s'}})}e$, where $({\bf{s}}, {\bf{s'}})$ denotes the number of pairings to the right and above of the unique pairing in ${\bf{s'}}$ but not in ${\bf{s}}$ (\cite[Page 1249, third paragraph in the proof of lemma 3.1]{bw} and \cite[Proof of Lemma 15]{w}).  For all $k$, $d_{\bf{n}}^{k}$ satisfies $d_{\bf{n}}^{k+1} \circ d_{\bf{n}}^{k}$ $=$ $0$ (See, \cite[Section 6.2]{w}).  

However, it has been unknown whether a Khovanov-type bicomplex exists for the homological grades $i$ and $k$ \cite[Section 3.3]{bw}.  
If such a bicomplex exists, there should be the spectral sequence whose $E_{2}$ term is determined by the bicomplex \cite[Section 3.3]{bw}.  

We now state the main results.  
\begin{mthm}\label{main-ex}
There exists a bicomplex $\{$$\mathcal{C}^{k, i}_{\bf{n}}(D)$, $d'^{k, i}$, $d''^{k, i}$$\}$.  
\end{mthm}
We prove the claim above in the following section.  
\section{Construction of the Khovanov-type bicomplex for the colored Jones polynomial.  }
In \cite{ito5}, we define the Khovanov homology $\mathcal{H}^{i}(D)$ $=$ $\mathcal{H}^{i}(\mathcal{C}^{*}(D), \delta_{s, t})$ of a link diagram $D$.  In this section, the differential $\delta_{s, t} : \mathcal{C}^{i}(D^{\bf{s}})$ $\to$ $\mathcal{C}^{i}(D^{\bf{s'}})$ is denoted by $d^{i}_{\bf{s}}$.  In order to get the bicomplex with the two homological degrees $k$ and $i$, let us define $d'^{k, i}_{\bf n} :$ $\mathcal{C}^{k, i}_{\bf{n}}(D)$ $\to$ $\mathcal{C}^{k+1, i}_{\bf{n}}(D)$.  Recall the map $e : {\bf{s}} \to {\bf{s'}}$ given by gluing an annulus between the two adjacent strands of the cable.  We will define the map $\mathcal{C}^{i}(D^{\bf{s}})$ $\to$ $\mathcal{C}^{i}(D^{\bf{s'}})$ corresponding to $e$.  In the following, we call the two adjacent strands that will be glued {\it contracted strands}.  We also call circles {\it contracted circles} if they consist of three or four contracted strands.  

\begin{figure}[h!]
\begin{minipage}{80pt}
\begin{picture}(70,80)
\put(17,0){(a-1)}
{\color{black}{
{\qbezier(23,54)(27,50)(31,46)}
}}
{\color{black}{
\qbezier(36,46)(40,50)(44,54)
}}
\put(0,50){\line(1,0){60}}
\put(20,20){\line(0,1){25}}
\put(40,20){\line(0,1){25}}
\put(20,55){\line(0,1){25}}
\put(40,55){\line(0,1){25}}
\put(30,75){\vector(0,-1){15}}
\end{picture}
\end{minipage}
\qquad
\begin{minipage}{80pt}
\begin{picture}(70,80)
{\color{black}{\qbezier(19,54)(23,50)(27,46)}}
{\color{black}{\qbezier(36,46)(40,50)(44,54)}}
\put(17,0){(a-2)}
\put(0,50){\line(1,0){15}}
\put(25,50){\line(1,0){10}}
\put(45,50){\line(1,0){15}}
\put(20,20){\line(0,1){60}}
\put(40,20){\line(0,1){60}}
\put(30,75){\vector(0,-1){15}}
\end{picture}
\end{minipage}
\qquad
\begin{minipage}{80pt}
\begin{picture}(70,80)
{\color{black}{\qbezier(27,64)(31,60)(35,56)}}
{\color{black}{\qbezier(24,36)(28,40)(32,44)}}
{\color{black}{\qbezier(40,56)(44,60)(48,64)}}
{\color{black}{\qbezier(36,44)(40,40)(44,36)}}
\put(20,0){(b)}
\put(0,40){\line(1,0){60}}
\put(0,60){\line(1,0){60}}
\put(20,20){\line(0,1){15}}
\put(40,20){\line(0,1){15}}
\put(20,45){\line(0,1){10}}
\put(40,45){\line(0,1){10}}
\put(20,65){\line(0,1){15}}
\put(40,65){\line(0,1){15}}
\end{picture}
\end{minipage}
\caption{(a-1), (a-2): Two crossings generated by two contracted strands and one non-contracted strand.  (b): Four crossings generated by only contracted strands.}\label{strand}
\end{figure}
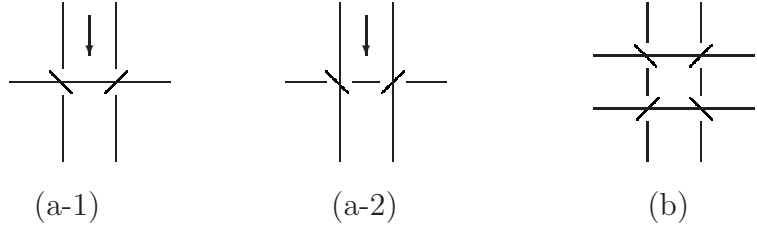

First, markers are put on $D^{\bf{s}}$ as in Figure \ref{strand}  
where the cases (a-1) and (a-2) depend on the orientations of contracted strands defined in the following.  
\begin{dfn}(the orientation of contracted strands)
Put $D^{\bf{s}}$ on $\mathbb{R}^{2}$ with a coordinate $(x_1, x_2)$ which has the only one maximum point for $x_2$-axis.  Let the base point be this maximum point and the orientation of the strands corresponding to the lower dot of $\Gamma_{\bf{n}}$ of the two.  Let $y$ be the word as the letters corresponding minus markers taken alternatively along the orientation of contracted strands from the base point.  
\end{dfn}
Note that either right or left crossing has a minus marker as in (a-1), (a-2) and (b) when we go along contracted strands and encounter another strands.  

\begin{figure}[h!]
\begin{minipage}{90pt}
\begin{picture}(80,90)
\put(42,46){\vector(0,-1){15}}
{\color{black}{\qbezier(28,14)(32,10)(36,6)}}
{\color{black}{\qbezier(24,74)(28,70)(32,66)}}
{\color{black}{\qbezier(40,6)(44,10)(48,14)}}
{\color{black}{\qbezier(36,66)(40,70)(44,74)}}
\put(-17,67){(1)}
\put(62,67){(2)}
\put(20,0){\line(0,1){60}}
\put(20,60){\line(0,1){7}}
\put(20,73){\line(0,1){7}}
\put(40,0){\line(0,1){60}}
\put(40,60){\line(0,1){7}}
\put(40,73){\line(0,1){7}}
\put(0,10){\line(1,0){15}}
\put(25,10){\line(1,0){10}}
\put(45,10){\line(1,0){15}}
\put(0,70){\line(1,0){60}}
\end{picture}
\end{minipage}
\caption{The figure shows a part of contracted strands going from crossings (a-1) to crossings (a-2) or (b).  All cases are got by considering every possible couple among (a-1), (a-2) and (b).  }\label{strand2}
\end{figure}
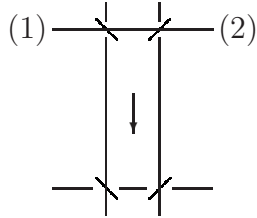

Second, we consider all possible choices of markers on crossings of non-contracted strands.  Also, we deal with disjoint circles by smoothing along these markers.  Put plus signs for all contracted circles and arbitrary signs for the other circles.  Automatically, an enhanced Kauffman state (so-called Type 1) is realised for each choice of markers and sings.  The contracted circles do not depend on the choice of markers of non-contracted circles.  This comes from the following fact.  Let us look at Figure \ref{strand2}.  
Strands (1) and (2) are edges of a circle by smoothening along the markers.  Then, a strand with (1) and (2) coincides with one of the enhanced Kauffman state $T_{1}$ of $D^{\bf{s'}}$ consisting of non-contracted circles by removing all of the contracted circles from Type 1.  Therefore contracted circles are unchanged up to plane isotopy for the choice of markers of the non-contracted strands.  

Finally, we show the existence of a differential of bicomplex.  Assume $(s, t)$ $\neq$ $(0, 0)$ for $\delta_{s, t}$.  We define Type 2 by using Type 1 in the following.  The markers of Type 2 is defined by replacing one pair 
\begin{minipage}{15pt}
\begin{picture}(15,10)
\put(0,5){\line(1,0){15}}
\put(11,0){\line(0,1){10}}
\put(4,0){\line(0,1){10}}
\qbezier(9,3)(11,5)(13,7)
\qbezier(2,7)(4,5)(6,3)
\end{picture}
\end{minipage} of Type 1 with 
\begin{minipage}{15pt}
\begin{picture}(15,10)
\put(0,5){\line(1,0){15}}
\put(11,0){\line(0,1){10}}
\put(4,0){\line(0,1){10}}
\qbezier(2,3)(4,5)(6,7)
\qbezier(9,7)(11,5)(13,3)
\end{picture}
\end{minipage}
where \begin{minipage}{15pt}
\begin{picture}(15,10)
\put(0,5){\line(1,0){15}}
\put(11,0){\line(0,1){10}}
\put(4,0){\line(0,1){10}}
\end{picture}
\end{minipage} stands for either 
\begin{minipage}{15pt}
\begin{picture}(15,10)
\put(0,5){\line(1,0){15}}
\put(11,0){\line(0,1){4}}
\put(11,10){\line(0,-1){4}}
\put(4,0){\line(0,1){4}}
\put(4,10){\line(0,-1){4}}
\end{picture}
\end{minipage} or 
\begin{minipage}{15pt}
\begin{picture}(15,10)
\put(0,5){\line(1,0){3}}
\put(5,5){\line(1,0){5}}
\put(12,5){\line(1,0){3}}
\put(11,0){\line(0,1){10}}
\put(4,0){\line(0,1){10}}
\end{picture}
\end{minipage}.  The correspondence of signs of Type 1 and Type 2 is obtained by Figure \ref{type-fig1} -- \ref{type1-2}.  Both (b-1) and (b-2) correspond to (a) in Figure \ref{type-fig1}.  In Figure \ref{type-fig2} and Figure \ref{type-fig3}, there are $1$--$1$ correspondences.  In Figure \ref{type-fig1} -- Figure \ref{type1-2}, $p:+$ (resp. $q:+$) stands for the multiplication of $p$ (resp. $q$) and $+$.  Note that $p:+$ and $q:+$ are not co-multiplications.  On the other hand, the pair $p:q$ and $q:p$ in Figure \ref{type-fig3} are Frobenius calculus defined in \cite[Section 1.2]{ito5}.  
\begin{figure}[!h]
\begin{minipage}{60pt}
\begin{picture}(60,100)
\qbezier(12,32)(17,37)(22,42)
\qbezier(12,68)(17,63)(22,58)
\qbezier(48,32)(43,37)(38,42)
\qbezier(38,58)(43,63)(48,68)
{\linethickness{0.4pt}
\put(-10,50){\vector(1,0){10}}
\put(30,95){\vector(0,-1){10}}}
\put(5,48){$p$}
\put(28,73){$q$}
\put(49,46){$+$}
\put(26,23){$+$}
\put(17,20){\line(0,1){60}}
\put(43,20){\line(0,1){60}}
\put(0,37){\line(1,0){15}}
\put(0,63){\line(1,0){15}}
\put(19,37){\line(1,0){22}}
\put(19,63){\line(1,0){22}}
\put(45,37){\line(1,0){15}}
\put(45,63){\line(1,0){15}}
\put(30,20){\oval(26,6)[b]}
\put(60,50){\oval(6,26)[r]}
\put(23,0){(a)}
\end{picture}
\end{minipage}
\qquad
\begin{minipage}{60pt}
\begin{picture}(60,100)
\qbezier(12,42)(17,37)(22,32)
\qbezier(12,68)(17,63)(22,58)
\qbezier(38,32)(43,37)(48,42)
\qbezier(38,58)(43,63)(48,68)
\put(-10,48){$p:+$}
\put(28,73){$q$}
\put(17,20){\line(0,1){60}}
\put(43,20){\line(0,1){60}}
\put(0,37){\line(1,0){15}}
\put(0,63){\line(1,0){15}}
\put(19,37){\line(1,0){22}}
\put(19,63){\line(1,0){22}}
\put(45,37){\line(1,0){15}}
\put(45,63){\line(1,0){15}}
\put(30,20){\oval(26,6)[b]}
\put(60,50){\oval(6,26)[r]}
\put(23,0){(b-1)}
\end{picture}
\end{minipage}
\qquad
\begin{minipage}{60pt}
\begin{picture}(60,100)
\qbezier(12,32)(17,37)(22,42)
\qbezier(12,68)(17,63)(22,58)
\qbezier(38,32)(43,37)(48,42)
\qbezier(38,68)(43,63)(48,58)
\put(5,48){$p$}
\put(18,73){$q:+$}
\put(17,20){\line(0,1){60}}
\put(43,20){\line(0,1){60}}
\put(0,37){\line(1,0){15}}
\put(0,63){\line(1,0){15}}
\put(19,37){\line(1,0){22}}
\put(19,63){\line(1,0){22}}
\put(45,37){\line(1,0){15}}
\put(45,63){\line(1,0){15}}
\put(30,20){\oval(26,6)[b]}
\put(60,50){\oval(6,26)[r]}
\put(23,0){(b-2)}
\end{picture}
\end{minipage}
\caption{(a): Type 1.  (b-1), (b-2): Type 2.  }\label{type-fig1}
\end{figure}
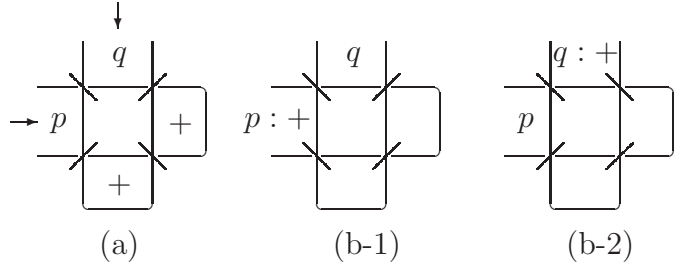
\begin{figure}[!h]
\begin{minipage}{60pt}
\begin{picture}(60,100)
\qbezier(12,32)(17,37)(22,42)
\qbezier(12,68)(17,63)(22,58)
\qbezier(48,32)(43,37)(38,42)
\qbezier(38,58)(43,63)(48,68)
\put(5,48){$p$}
\put(28,73){$q$}
\put(49,46){$+$}
\put(26,23){$+$}
\put(0,37){\line(1,0){60}}
\put(0,63){\line(1,0){60}}
\put(17,20){\line(0,1){15}}
\put(43,20){\line(0,1){15}}
\put(17,39){\line(0,1){22}}
\put(43,39){\line(0,1){22}}
\put(17,65){\line(0,1){15}}
\put(43,65){\line(0,1){15}}
\put(30,20){\oval(26,6)[b]}
\put(60,50){\oval(6,26)[r]}
\put(23,0){(a)}
{\linethickness{0.4pt}
\put(-10,50){\vector(1,0){10}}
\put(30,95){\vector(0,-1){10}}}
\end{picture}
\end{minipage}
\qquad
\begin{minipage}{60pt}
\begin{picture}(60,100)
\qbezier(12,32)(17,37)(22,42)
\qbezier(12,68)(17,63)(22,58)
\qbezier(38,32)(43,37)(48,42)
\qbezier(38,68)(43,63)(48,58)
\put(5,48){$p$}
\put(18,73){$q:+$}
\put(0,37){\line(1,0){60}}
\put(0,63){\line(1,0){60}}
\put(17,20){\line(0,1){15}}
\put(43,20){\line(0,1){15}}
\put(17,39){\line(0,1){22}}
\put(43,39){\line(0,1){22}}
\put(17,65){\line(0,1){15}}
\put(43,65){\line(0,1){15}}
\put(30,20){\oval(26,6)[b]}
\put(60,50){\oval(6,26)[r]}
\put(23,0){(b)}
\end{picture}
\end{minipage}
\caption{(a): Type 1.  (b): Type 2.  }\label{type-fig2}
\end{figure}
\begin{figure}[!h]
\begin{minipage}{60pt}
\begin{picture}(60,100)
\qbezier(12,32)(17,37)(22,42)
\qbezier(12,68)(17,63)(22,58)
\qbezier(48,32)(43,37)(38,42)
\qbezier(38,58)(43,63)(48,68)
\put(-5,48){$p:q$}
\put(18,73){$q:p$}
\put(49,46){$+$}
\put(26,23){$+$}
\put(0,37){\line(1,0){60}}
\put(0,63){\line(1,0){60}}
\put(17,20){\line(0,1){15}}
\put(43,20){\line(0,1){15}}
\put(17,39){\line(0,1){22}}
\put(43,39){\line(0,1){22}}
\put(17,65){\line(0,1){15}}
\put(43,65){\line(0,1){15}}
\put(30,20){\oval(26,6)[b]}
\put(60,50){\oval(6,26)[r]}
\put(23,0){(a)}
{\linethickness{0.4pt}
\put(-20,50){\vector(1,0){10}}
\put(30,95){\vector(0,-1){10}}}
\end{picture}
\end{minipage}
\qquad
\begin{minipage}{60pt}
\begin{picture}(60,100)
\qbezier(22,32)(17,37)(12,42)
\qbezier(12,58)(17,63)(22,68)
\qbezier(48,32)(43,37)(38,42)
\qbezier(38,58)(43,63)(48,68)
\put(6,70){$p$}
\put(18,73){$q:+$}
\put(49,46){$+$}
\put(0,37){\line(1,0){60}}
\put(0,63){\line(1,0){60}}
\put(17,20){\line(0,1){15}}
\put(43,20){\line(0,1){15}}
\put(17,39){\line(0,1){22}}
\put(43,39){\line(0,1){22}}
\put(17,65){\line(0,1){15}}
\put(43,65){\line(0,1){15}}
\put(30,20){\oval(26,6)[b]}
\put(60,50){\oval(6,26)[r]}
\put(23,0){(b)}
\end{picture}
\end{minipage}
\caption{(a): Type 1.  (b): Type 2.  }\label{type-fig3}
\end{figure}
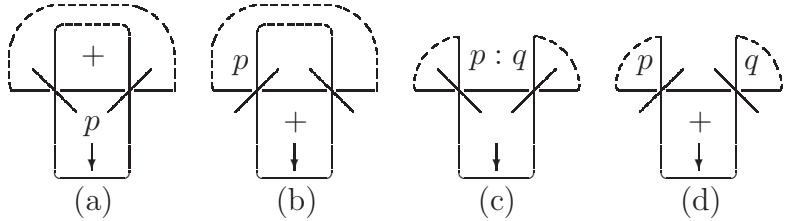
\begin{figure}[h!]
\begin{minipage}{60pt}
\begin{picture}(60,80)
\put(0,40){\line(1,0){14}}
\put(18,40){\line(1,0){24}}
\put(46,40){\line(1,0){14}}
\put(44,20){\line(0,1){40}}
\put(16,20){\line(0,1){40}}
\qbezier(36,32)(44,40)(52,48)
\qbezier(8,48)(16,40)(24,32)
\put(26,50){$+$}
\put(27,25){$p$}
\qbezier(16,19)(16,18)(16,17)
\qbezier(16,16)(16,15)(16,14)
\qbezier(16,13)(16,12)(16,11)
\qbezier(44,19)(44,18)(44,17)
\qbezier(44,16)(44,15)(44,14)
\qbezier(44,13)(44,12)(44,11)
\put(30,10){\oval(28,6)[b]}
\qbezier(16,60.8)(16,61.3)(16,62)
\qbezier(16,62)(16,62.5)(16.5,63.2)
\qbezier(17,64)(18,65)(19,65)
\qbezier(20.3,65)(21.3,65)(22.3,65)
\qbezier(23.8,65)(24.8,65)(25.8,65)
\qbezier(27.6,65)(28.6,65)(29.6,65)
\qbezier(31,65)(32,65)(33,65)
\qbezier(34.3,65)(35.3,65)(36.3,65)
\qbezier(38,65)(39,65)(40,65)
\qbezier(41.5,64.8)(42,64.8)(43,64)
\qbezier(43.6,63)(43.8,63.2)(44,61)
\put(23,-5){(a)}
\qbezier(0,40)(-0.5,40)(-1,43)
\qbezier(-1,45)(-1,46)(-1,47)
\qbezier(-1,49)(-1,50)(-1,51)
\qbezier(-1,53)(-1,54)(-1,55)
\qbezier(-0.7,57)(0,59)(0.8,61)
\qbezier(1.8,63)(2.3,64)(3.3,65)
\qbezier(5,67)(6,68)(7.5,69)
\qbezier(9.5,70.2)(10.5,70.8)(11.5,71.1)
\qbezier(13.5,71.8)(14.5,72)(15.5,72.5)
\qbezier(17.5,72.5)(18.5,72.5)(19.5,72.5)
\qbezier(21.5,72.5)(22.5,72.5)(23.5,72.5)
\qbezier(25.5,72.5)(26.5,72.5)(27.5,72.5)
\qbezier(29.5,72.5)(30.5,72.5)(31.5,72.5)
\qbezier(33.5,72.5)(34.5,72.5)(35.5,72.5)
\qbezier(37.5,72.5)(38.5,72.5)(39.5,72.5)
\qbezier(41.5,72.5)(42.5,72.5)(43.5,72.5)
\qbezier(45.5,72)(46.5,71.8)(47.5,71.4)
\qbezier(48.5,71)(49.5,70.6)(50.5,70.1)
\qbezier(52.5,68.8)(53.5,68.3)(54.5,67.2)
\qbezier(56.5,65.2)(57.5,64)(58.5,62.3)
\qbezier(59.4,60)(60,59)(60.2,58)
\qbezier(60.6,56)(60.9,55)(60.8,54)
\qbezier(61,52)(61,51)(61,50)
\qbezier(61,48)(61,46.5)(61,45)
\qbezier(61,43)(61,40)(61,40)
{\linethickness{0.4pt}
\put(30,20){\vector(0,-1){10}}}
\end{picture}
\end{minipage}
\quad
\begin{minipage}{60pt}
\begin{picture}(60,80)
\put(0,40){\line(1,0){14}}
\put(18,40){\line(1,0){24}}
\put(46,40){\line(1,0){14}}
\put(44,20){\line(0,1){40}}
\put(16,20){\line(0,1){40}}
\qbezier(36,48)(44,40)(52,32)
\qbezier(8,32)(16,40)(24,48)
\put(7,48){$p$}
\put(26,25){$+$}
\qbezier(16,19)(16,18)(16,17)
\qbezier(16,16)(16,15)(16,14)
\qbezier(16,13)(16,12)(16,11)
\qbezier(44,19)(44,18)(44,17)
\qbezier(44,16)(44,15)(44,14)
\qbezier(44,13)(44,12)(44,11)
\put(30,10){\oval(28,6)[b]}
\qbezier(16,60.8)(16,61.3)(16,62)
\qbezier(16,62)(16,62.5)(16.5,63.2)
\qbezier(17,64)(18,65)(19,65)
\qbezier(20.3,65)(21.3,65)(22.3,65)
\qbezier(23.8,65)(24.8,65)(25.8,65)
\qbezier(27.6,65)(28.6,65)(29.6,65)
\qbezier(31,65)(32,65)(33,65)
\qbezier(34.3,65)(35.3,65)(36.3,65)
\qbezier(38,65)(39,65)(40,65)
\qbezier(41.5,64.8)(42,64.8)(43,64)
\qbezier(43.6,63)(43.8,63.2)(44,61)
\put(23,-5){(b)}
\qbezier(0,40)(-0.5,40)(-1,43)
\qbezier(-1,45)(-1,46)(-1,47)
\qbezier(-1,49)(-1,50)(-1,51)
\qbezier(-1,53)(-1,54)(-1,55)
\qbezier(-0.7,57)(0,59)(0.8,61)
\qbezier(1.8,63)(2.3,64)(3.3,65)
\qbezier(5,67)(6,68)(7.5,69)
\qbezier(9.5,70.2)(10.5,70.8)(11.5,71.1)
\qbezier(13.5,71.8)(14.5,72)(15.5,72.5)
\qbezier(17.5,72.5)(18.5,72.5)(19.5,72.5)
\qbezier(21.5,72.5)(22.5,72.5)(23.5,72.5)
\qbezier(25.5,72.5)(26.5,72.5)(27.5,72.5)
\qbezier(29.5,72.5)(30.5,72.5)(31.5,72.5)
\qbezier(33.5,72.5)(34.5,72.5)(35.5,72.5)
\qbezier(37.5,72.5)(38.5,72.5)(39.5,72.5)
\qbezier(41.5,72.5)(42.5,72.5)(43.5,72.5)
\qbezier(45.5,72)(46.5,71.8)(47.5,71.4)
\qbezier(48.5,71)(49.5,70.6)(50.5,70.1)
\qbezier(52.5,68.8)(53.5,68.3)(54.5,67.2)
\qbezier(56.5,65.2)(57.5,64)(58.5,62.3)
\qbezier(59.4,60)(60,59)(60.2,58)
\qbezier(60.6,56)(60.9,55)(60.8,54)
\qbezier(61,52)(61,51)(61,50)
\qbezier(61,48)(61,46.5)(61,45)
\qbezier(61,43)(61,40)(61,40)
{\linethickness{0.4pt}
\put(30,20){\vector(0,-1){10}}}
\end{picture}
\end{minipage}
\quad
\begin{minipage}{60pt}
\begin{picture}(60,80)
\put(0,40){\line(1,0){14}}
\put(18,40){\line(1,0){24}}
\put(46,40){\line(1,0){14}}
\put(44,20){\line(0,1){40}}
\put(16,20){\line(0,1){40}}
\qbezier(36,32)(44,40)(52,48)
\qbezier(8,48)(16,40)(24,32)
\put(20,50){$p:q$}
\qbezier(16,19)(16,18)(16,17)
\qbezier(16,16)(16,15)(16,14)
\qbezier(16,13)(16,12)(16,11)
\qbezier(44,19)(44,18)(44,17)
\qbezier(44,16)(44,15)(44,14)
\qbezier(44,13)(44,12)(44,11)
\put(30,10){\oval(28,6)[b]}
\qbezier(0.5,40)(-0.8,40)(-1,42)
\qbezier(-0.8,44)(-0.7,45)(-0.3,46)
\qbezier(0.2,48)(0.5,49)(1.2,50)
\qbezier(2.2,52)(2.8,53)(3.8,54)
\qbezier(6,56.1)(7,57.2)(8,57.6)
\qbezier(10,58.8)(11,59.3)(12,59.8)
\qbezier(14,60.3)(16,60.5)(16,59.5)
\qbezier(44,59)(44,61)(45,60.5)
\qbezier(46,60.2)(47,60)(48,59.6)
\qbezier(50,58.8)(51,58.5)(52,57.6)
\qbezier(54,56.1)(55,55.5)(56.5,53.6)
\qbezier(58,51.5)(59,50)(59.5,48.3)
\qbezier(60.4,46)(60.7,45)(60.8,44)
\qbezier(61,42)(61,40)(59.8,40)
\put(23,-5){(c)}
{\linethickness{0.4pt}
\put(30,20){\vector(0,-1){10}}}
\end{picture}
\end{minipage}
\quad
\begin{minipage}{60pt}
\begin{picture}(60,80)
\put(0,40){\line(1,0){14}}
\put(18,40){\line(1,0){24}}
\put(46,40){\line(1,0){14}}
\put(44,20){\line(0,1){40}}
\put(16,20){\line(0,1){40}}
\qbezier(36,48)(44,40)(52,32)
\qbezier(8,32)(16,40)(24,48)
\put(7,48){$p$}
\put(47,48){$q$}
\put(26,25){$+$}
\qbezier(16,19)(16,18)(16,17)
\qbezier(16,16)(16,15)(16,14)
\qbezier(16,13)(16,12)(16,11)
\qbezier(44,19)(44,18)(44,17)
\qbezier(44,16)(44,15)(44,14)
\qbezier(44,13)(44,12)(44,11)
\put(30,10){\oval(28,6)[b]}
\qbezier(0.5,40)(-0.8,40)(-1,42)
\qbezier(-0.8,44)(-0.7,45)(-0.3,46)
\qbezier(0.2,48)(0.5,49)(1.2,50)
\qbezier(2.2,52)(2.8,53)(3.8,54)
\qbezier(6,56.1)(7,57.2)(8,57.6)
\qbezier(10,58.8)(11,59.3)(12,59.8)
\qbezier(14,60.3)(16,60.5)(16,59.5)
\qbezier(44,59)(44,61)(45,60.5)
\qbezier(46,60.2)(47,60)(48,59.6)
\qbezier(50,58.8)(51,58.5)(52,57.6)
\qbezier(54,56.1)(55,55.5)(56.5,53.6)
\qbezier(58,51.5)(59,50)(59.5,48.3)
\qbezier(60.4,46)(60.7,45)(60.8,44)
\qbezier(61,42)(61,40)(59.8,40)
\put(23,-5){(d)}
{\linethickness{0.4pt}
\put(30,20){\vector(0,-1){10}}}
\end{picture}
\end{minipage}
\caption{(a)--(d) are concerned with contracted strands.  Arrows stand for orientations of contracted strands.  In these figure, $p$ and $q$ stand for signs of Kauffman states.  (a): A part of Type 1 corresponding to (b).  (b): A part of Type 2 generating two circles.  (c): A part of Type 1 corresponding to (d).  (d): A part of Type 2 generating three circles.  The symbol $p:q$ is the multiplication of two circles with $p$ and $q$.  }\label{type1-2}
\end{figure}
We define $d^{k, i}_{\bf n} :$ $\mathcal{C}^{k, i}_{\bf n}(D)$ $\to$ $\mathcal{C}^{k+1, i}_{\bf n}(D)$ by $\mathcal{C}^{i}(D^{\bf s}) \ni$ $S \otimes [xy]$ $\mapsto$ $S' \otimes [x]$ $\in$ $\mathcal{C}^{i}(D^{\bf s'})$ where $S'$ is $T_{1}$ (resp. $- T_{1}$) if $S$ is Type 1 (resp. Type 2) 0 otherwise.  By the definition, we have $d_{\bf{s'}}^{i} \circ d_{\bf{n}}^{k, i}$ $=$ $d_{\bf{n}}^{k, i+1} \circ d_{\bf{s}}^{i}$.  
Setting $d'^{k, i}$ $:=$ $(-1)^{({\bf{s}}, {\bf{s'}})}d_{\bf{n}}^{k, i}$ and $d''^{k, i}$ $:=$ $(-1)^{k}$ $\oplus_{{\bf{s}} \in I_{\bf{k}}, |{\bf{k}}| = k}d_{\bf{s}}^{i}$, we have $d''^{k+1, i} \circ d'^{k, i}$ $+$ $d'^{k, i+1} \circ d''^{k, i}$ $=$ $0$.  This completes the proof of Theorem \ref{main-ex}.  

\begin{rem}
The explicit chain homotopy maps implying the second and third Reidemeister invariance of the Khovanov homology of $\mathcal{C}_{\bf{n}}^{k, i, j}(D)$ with the differential defined in \cite[Section 2.2]{jacobsson} are given by \cite[Equation (7), (12)]{ito3}.  
\end{rem}

\begin{rem}
By checking all types of crossings, we notice that the number of contracted circles are even.  Instead of putting minus signs for all contracted strands in the proof above, we consider every enhanced state which has equal numbers of positive contracted circles and minus circles and set $\delta_{s, t}$ $=$ $\delta_{0, 0}$.  In this case, the following fact is available where $d'^{k, i, j}$ $:=$ $d'^{k, i}|_{\mathcal{C}^{k, i, j}_{\bf n}(D)}$ and $d''^{k, i, j}$ $:=$ $d''^{k, i}|_{\mathcal{C}^{k, i, j}_{\bf n}(D)}$: 
\begin{thm}
There exists a bicomplex $\{\mathcal{C}^{k, i, j}_{\bf n}(D), d'^{k, i, j}, d''^{k, i, j}\}$.  
\end{thm}
\end{rem}

\scriptsize{\textsc{Department of Mathematics Waseda University.  Tokyo 169-8555, Japan.  }}
\scriptsize{{\it E-mail address:} \texttt{noboru@moegi.waseda.jp}}
\end{document}